\documentclass[11pt]{article}
\usepackage{amsmath,amssymb, amsbsy}
\usepackage{color,psfrag}
\usepackage[dvips]{graphicx}
\usepackage{enumerate}
\newcommand{\R}{{\mathbb R}}
\newcommand{\N}{{\mathbb N}}

\newcommand{\eps}{\varepsilon}

\newcommand{\E}{{\mathcal E}}
\renewcommand{\phi}{\varphi}
\def\endproof{\hfill $\Box$\par\vskip3mm}

\renewcommand{\ge }{\geqslant}
\renewcommand{\geq }{\geqslant}
\renewcommand{\le }{\leqslant}
\renewcommand{\leq }{\leqslant}
\def\neweq#1{\begin{equation}\label{#1}}
\def\endeq{\end{equation}}
\def\eq#1{(\ref{#1})}
\def\proof{{\em Proof.} }

\DeclareMathOperator{\cn}{cn}
\DeclareMathOperator{\sn}{sn}
\DeclareMathOperator{\dn}{dn}

\newtheorem{theorem}{Theorem}
\newtheorem{proposition}[theorem]{Proposition}

\newtheorem{corollary}[theorem]{Corollary}

\newtheorem{definition}{Definition}

\begin{document}

\title{Resonance tongues for the Hill equation with Duffing coefficients\\
and instabilities in a nonlinear beam equation}

\author{Carlo GASPARETTO -- Filippo GAZZOLA\\
{\small Dipartimento di Matematica del Politecnico, Piazza L. da Vinci 32 - 20133 Milano (Italy)}}

\date{}
\maketitle

\begin{abstract}
We consider a class of Hill equations where the periodic coefficient is the squared solution of some Duffing equation plus a constant.
We study the stability of the trivial solution of this Hill equation and we show that a criterion due to Burdina \cite{burdina} is very
helpful for this analysis. In some cases, we are also able to determine exact solutions in terms of Jacobi elliptic functions.
Overall, we obtain a fairly complete picture of the stability and instability regions. These results are then used to study the
stability of nonlinear modes in some beam equations.
\end{abstract}

\section{Introduction}

In order to explain the contents of the present paper, we briefly introduce the Duffing and Hill equations.
The Duffing equation \cite{duffing} (see also \cite{stoker}) is a nonlinear ODE and reads
\neweq{ODE}
\ddot{y}(t)+y(t)+y(t)^3=0\qquad(t>0)\, .
\endeq
To \eq{ODE} we associate the initial values
\neweq{alphabeta}
y(0)=\delta\ ,\qquad \dot{y}(0)=0\, ,\qquad(\delta\in\R\setminus\{0\})\, .
\endeq
The unique solution of \eq{ODE}-\eq{alphabeta} is periodic, its period depends on $\delta$ and may be computed in terms of Jacobi elliptic functions,
see Proposition \ref{prop:sol_duffing} in Section \ref{properties}.\par
The Hill equation \cite{Hill} was introduced for the study of the lunar perigee and it has been the object of many subsequent studies,
see e.g.\ \cite{cesari,magnus,stoker,yakubovich}. It is a linear ODE with periodic coefficient, that is,
\neweq{hill}
\ddot{\xi}(t)+p(t)\xi(t)=0\, ,\qquad p\in C^0[0,T]\, ,\quad p(t+T)=p(t)\quad \forall t
\endeq
where we intend that $T>0$ is the smallest period of $p$ (in particular, $p$ is nonconstant). The main concern is to establish whether the trivial
solution $\xi\equiv0$ of \eq{hill} is stable or, equivalently, if all the solutions of \eq{hill} are bounded in $\R$. If $p(t)=a+2q\cos(2t)$
for some $a,q>0$, then \eq{hill} is named after Mathieu \cite{mathieu} and, in this case, the stability analysis for \eq{hill} is well-understood:
in the $(q,a)$-plane, the so-called {\em resonance tongues} (or instability regions) for \eq{hill} emanate from the points $(0,\ell^2)$, with $\ell\in\N$,
see \cite[fig.8A]{mcl}: these tongues are separated from the stability regions by some {\em resonance lines} which are explicitly known.
This is one of the few cases where the stability for \eq{hill} has reached a complete understanding. In general, the resonance tongues
have strange geometries, see e.g.\ \cite{broer,broer2,simakhina,svetlana,weinstein}.\par
The first purpose of the present paper is to study the stability of the Hill equation \eq{hill} when the periodic coefficient $p$ is related to
a solution of the Duffing equation \eq{ODE}. More precisely, we will consider the cases where
\neweq{where}
p(t)=\gamma+\beta y(t)^2
\endeq
with $\gamma\ge0$, $\beta>0$, and $y$ being a solution of \eq{ODE} or of a scaled version of it. We first consider the simpler case
where $\beta=1$. Since the solution of \eq{ODE} is even with respect to $t$, the stability diagram for \eq{hill} in the $(\delta,\gamma)$-plane is
symmetric with respect to the axis $\delta=0$. By allowing $\gamma$ to become negative we obtain the explicit form of the first
resonance tongue and the stability behavior of \eq{hill} in a large neighborhood of $(\delta,\gamma)=(0,0)$, see Theorems \ref{exactsolutions}
and \ref{striscia}, as well as Figure \ref{zonevere}. In our analysis we take advantage of some explicit solutions of \eq{hill} that we
are able to determine thanks to the particular form of the coefficient $p$ in \eq{where} when $\beta=1$.\par
In Theorem \ref{primo} we find sufficient conditions
on $\gamma$ and $\delta$ for the stability of \eq{hill} when $p$ is as in \eq{where} and $\beta=1$. The proof contains two main ideas. First, we apply the
Burdina criterion \cite{burdina}, reported in Proposition \ref{lyapzhu} $(iii)$. Second, we prove that this criterion {\em does not} need the explicit
form of the solution $y$ of \eq{ODE} and the stability analysis may be reduced to the study of some elliptic integrals. In order to show that the
sufficient conditions are accurate, we proceed numerically. In Figure \ref{num_B} we plot the stability regions obtained through the Burdina criterion
whereas in Figure \ref{realplot} we plot the stability regions
obtained through the trace of the monodromy matrix. These plots confirm that the Burdina criterion is indeed quite accurate. Probably, it is
the most accurate for the Hill equation with squared Duffing coefficients: in Figure \ref{num_B} we also compare the Burdina criterion with the
Zhukovskii criterion \cite{zhk} (reported in Proposition \ref{lyapzhu} $(ii)$) and we show that the former is much more precise.\par
The Hill equation \eq{hill} with squared Duffing coefficients \eq{where} for general $\beta>0$ arises when studying the stability of modes
of nonlinear strings \cite{cazw,cazw2,ghg}, beams \cite{bbfg}, and plates \cite{bfg1,fergazmor}. In 1950, Woinowsky-Krieger
\cite{woinowsky} modified the classical beam models by Bernoulli and Euler assuming a nonlinear dependence of the axial strain on the
deformation gradient, by taking into account the stretching of the beam due to its elongation. Independently, Burgreen \cite{burg}
derived the very same nonlinear beam equation. After normalization of some physical constants and after scaling the space variable, the problem reads
\neweq{truebeam}
\left\{\begin{array}{ll}
u_{tt}+u_{xxxx}-\frac{2}{\pi}\, \|u_x\|^2_{L^2(0,\pi)}\, u_{xx}=0\quad & x\in(0,\pi)\, ,\ t>0\, ,\\
u(0,t)=u(\pi,t)=u_{xx}(0,t)=u_{xx}(\pi,t)=0\quad & t>0\, ,
\end{array}\right.
\endeq
where the Navier boundary conditions model a beam hinged at its endpoints. The term $\frac{2}{\pi}\|u_x\|^2_{L^2(0,\pi)}$ measures
the geometric nonlinearity of the beam due to its stretching.
In Section \ref{appl}, we recall the definitions of nonlinear modes of \eq{truebeam} and of their linear stability
and how this study naturally leads to \eq{hill} with $p$ as in \eq{where} with $\beta\neq1$.
The linear stability for problem \eq{truebeam} was recently tackled in \cite{bbfg} by using both the Zhukovskii criterion \cite{zhk}
and a generalization of the Lyapunov \cite{lyapunov} criterion due to Li-Zhang \cite[Theorem 1]{lizhang}. These criteria did not allow a full
understanding of the linear stability for \eq{truebeam} and several problems had to be left open. In this paper we apply
the Burdina criterion \cite{burdina} and we show that it allows to solve some of the open problems and to give a better description of the resonance tongues
for couples of modes of \eq{truebeam}.\par
This paper is organized as follows. In the next section we recall some basic facts about equations \eq{ODE} and \eq{hill}. In Section \ref{duffsqMOD}
we state our results about the Hill equation $\ddot{\xi}(t)+\big(\gamma+y(t)^2\big)\xi(t)=0$, where $y$ is the solution of \eq{ODE}-\eq{alphabeta}.
In Section \ref{appl} we define the nonlinear modes of \eq{truebeam} and we perform the same analysis for a more general Hill equation in order to study
the stability of these modes. All the proofs are postponed until Section~\ref{proofs}. In the final section we list some open problems.

\section{Basic properties of the Duffing and the Hill equations}\label{properties}

The Duffing equation \eq{ODE} admits periodic solutions whose period depends on the parameter $\delta$ in \eq{alphabeta}. An explicit solution of~\eqref{ODE},
which uses elliptic Jacobi functions, has been given by Burgreen~\cite{burg}: its period may be computed by using some properties of elliptic functions and
without knowing the explicit form of it. We refer to \cite{abram} for the definitions and the basic properties of the elliptic functions and integrals.
For the sake of completeness and because we need a formula therein, we briefly sketch a different proof of the computation of the period.

\begin{proposition} \label{prop:sol_duffing}
For all $\delta\in\R$ the solution of the Cauchy problem~\eqref{ODE}-\eqref{alphabeta} is
\begin{equation}\label{soluzione_duffing}
y(t) = \delta \cn{\bigg[t\sqrt{1+\delta^2},\frac{\delta}{\sqrt{2(1+\delta^2)}}\bigg]}\, ,
\end{equation}
where $\cn{[u,k]}$ is the Jacobi elliptic cosine function. Moreover, the period of the solution~\eqref{soluzione_duffing} is given by
\neweq{TE}
T(\delta)=4\sqrt2\int_0^1\frac{d\theta}{\sqrt{(2+\delta^2+\delta^2\theta^2)(1-\theta^2)}}\, .
\endeq
\end{proposition}
\proof The equation \eq{ODE} is conservative, its solutions have constant energy which is defined by
\neweq{energyy}
E(\delta)=\frac{\dot{y}^2}{2}+\frac{y^2}{2}+\frac{y^4}{4}\equiv\frac{\delta^2}{2}+\frac{\delta^4}{4}\, ,
\endeq
where $\delta$ is the initial condition in \eq{alphabeta}. We rewrite \eq{energyy} as
\neweq{ancoraenergia}
2\dot{y}^2=(2+\delta^2+y^2)(\delta^2-y^2)\qquad\forall t
\endeq
so that $-|\delta|\le y(t)\le|\delta|$ for all $t$
and $y$ oscillates in this range. If $y(t)$ solves \eq{ODE}-\eq{alphabeta}, then also $y(-t)$ solves the same problem:
this shows that the period $T(\delta)$ of $y$ is the double of the length of an interval of monotonicity for $y$.
Since the problem is autonomous, we may assume that $y(0)=-\delta<0$ and $\dot{y}(0)=0$; then we have that
$y(T/2)=\delta$ and $\dot{y}(T/2)=0$. By rewriting \eq{ancoraenergia} as
\neweq{derivata}
\sqrt2 \dot{y}=\sqrt{(2+\delta^2+y^2)(\delta^2-y^2)}\qquad\forall t\in\left(0,\frac{T}{2}\right)\, ,
\endeq
by separating variables, and upon integration over the time interval $(0,T/2)$ we obtain
$$\frac{T(\delta)}{2}=\sqrt2 \int_{-\delta}^{\delta}\frac{d\, y}{\sqrt{(2+\delta^2+y^2)(\delta^2-y^2)}}\,.$$
Then, using the fact that the integrand is even with respect to $y$ and through the change of variable $y=\delta\theta$, we obtain \eq{TE}.\endproof

Proposition \ref{prop:sol_duffing} tells that the period of $y$ depends on its initial value $\delta$ and therefore on its energy \eq{energyy}: this
is clearly due to the nonlinear nature of \eq{ODE}.
We point out that a solution $y$ of \eq{ODE} has mean value 0 over a period and that $y(t)$ and its translation $y(t+t_0)$ all have the same
energy \eq{energyy} for any $t_0\in\R$; therefore, initial conditions different from \eq{alphabeta}, but with the same energy, lead to the same
solution of \eq{ODE}, up to a time translation.\par
From \eq{TE} we infer that the map $\delta\mapsto T(\delta)$ is strictly decreasing on $(0,+\infty)$ and
$$
\lim_{\delta\downarrow0}T(\delta)=2\pi\, ,\qquad\lim_{\delta\to+\infty}T(\delta)=0\, .
$$
The period $T(\delta)$ can also be expressed differently. With the change of variables $\theta=\cos\alpha$,
\eq{TE} becomes
\neweq{firstkind}
T(\delta)=4\sqrt2 \int_0^{\pi/2}\frac{d\alpha}{\sqrt{2+\delta^2+\delta^2\cos^2\alpha}}=
\frac{4}{\sqrt{1+\delta^2}}\, K\left(\frac{\delta}{\sqrt{2(1+\delta^2)}}\right)\, ,
\endeq
where $K(\cdot)$ is the complete elliptic integral of the first kind, see \cite{abram}.\par
A further constant which plays an important role in our analysis is
\begin{equation}\label{costantesigma}
\sigma:=\int_0^1\frac{d\theta}{\sqrt{1-\theta^4}}=\int_0^{\pi/2}\frac{d\alpha}{\sqrt{1+\sin^2{\alpha}}}=\frac{1}{\sqrt{2}}K\Bigg(\frac{1}{\sqrt{2}}\Bigg)
=\frac{\sqrt{\pi}}{4}\frac{\Gamma(\frac{1}{4})}{\Gamma(\frac{3}{4})}\, .
\end{equation}

It is in general extremely difficult to establish whether the parameters of the Hill equation \eq{hill} lead to a stable regime.
In most cases, one uses suitable criteria which yield sufficient conditions for stability. There are many such criteria, see
e.g.\ \cite{cesari,magnus,stoker} and references therein. We state here three criteria which appear appropriate to tackle our problem.

\begin{proposition}\label{lyapzhu}
Assume that one of the three following facts holds:\par
$$(i)\ p\ge0\ \mbox{ and }\ T^3\int_0^T p(t)^2\, dt< \frac{64}{3}\sigma^4\qquad(\sigma\mbox{ as in \eqref{costantesigma}})\, ,$$
$$(ii)\ p\ge0\mbox{ and }
\exists\ell\in\N\mbox{ s.t. }\frac{\ell^2\pi^2}{T^2}\le p(t)\le\frac{(\ell+1)^2\pi^2}{T^2}\quad\forall t\, ,$$
$$(iii)\ p>0,\ p\mbox{ admits a unique maximum point and a unique minimum point in }[0,T),\mbox{ and }$$
$$\exists\ell\in\N\mbox{ s.t. }\ \ell\pi<\int_0^T\!\!\sqrt{p(t)}\, dt-\tfrac12\log\tfrac{\max p}{\min p}\le
\int_0^T\!\!\sqrt{p(t)}\, dt+\tfrac12\log\tfrac{\max p}{\min p}<(\ell+1)\pi\, .$$
Then the trivial solution of \eqref{hill} is stable.
\end{proposition}
\proof The first criterion is due to Li-Zhang \cite{lizhang} and generalizes the original Lyapunov criterion \cite{lyapunov}.
The second criterion is due to Zhukovskii \cite{zhk}; see also \cite[Test 1, \S 3, Chapter VIII]{yakubovich}.
The third criterion is due to Burdina \cite{burdina}; see also \cite[Test 3, \S 3, Chapter VIII]{yakubovich}.\endproof

\section{Stability for the Hill equation with a squared Duffing coefficient}\label{duffsqMOD}

In this section, we analyze the stability regions for the equation

\begin{equation}\label{squareduffing}
\ddot{\xi}(t)+\Big(\gamma+y(t)^2\Big)\xi(t)=0
\end{equation}
in the parameter $(\delta,\gamma)$-plane; here $y$ is the solution of~\eqref{ODE}-\eqref{alphabeta}. As a straightforward consequence of
Proposition \ref{prop:sol_duffing}, we infer that if $y$ solves \eqref{ODE}-\eqref{alphabeta}, then the squared function $y^2$
is periodic and its period is given by $T(\delta)/2$, that is,
$$
\frac{T(\delta)}{2}=2\sqrt2\int_0^1\frac{d\theta}{\sqrt{(2+\delta^2+\delta^2\theta^2)(1-\theta^2)}}\, .
$$

Our first result concerns the limit case $\gamma=0$.
\begin{theorem}\label{intornorigine}
	If $\gamma=0$, then the trivial solution of~\eqref{squareduffing} is stable for all $\delta>0$.
\end{theorem}

This result is a particular case of Theorem~\ref{striscia} below. However, since our proof of Theorem~\ref{intornorigine} makes use of the Li-Zhang \cite{lizhang} generalized Lyapunov criterion, it has its own independent interest and we give its proof in Section~\ref{proofs}.

As we know from~\cite[Chapter VIII]{yakubovich}, for all $\ell\in\N$, there exists a resonant tongue $U_\ell$ emanating from the point $(\delta,\gamma)=(0,\ell^2)$. If $(\delta,\gamma)$ belongs to one of these tongues, then the trivial solution of \eqref{squareduffing} is unstable.
In Section~\ref{proofs}, we prove the following precise characterization of the first instability region $U_1$.

\begin{theorem}\label{exactsolutions}
Let $U_1$ be the resonant tongue of \eqref{squareduffing} emanating from $(\delta,\gamma)=(0,1)$. Then
$$U_1=\left\{(\delta,\gamma)\in\R^2_+;\, 1<\gamma<1+\frac{\delta^2}{2}\right\}\, .$$
\end{theorem}

In particular, Theorem \ref{intornorigine} states that, contrary to what happens for the Mathieu equations, the resonance lines do not bend
downwards since they remain above the parabola $\gamma=1+\delta^2/2$. Moreover, as a consequence of Theorems~\ref{intornorigine}
and~\ref{exactsolutions}, we see that the strip $0\leq\gamma<1$ is a region of stability for~\eqref{squareduffing}; this follows by applying Theorem II,
p.695 in~\cite{yakubovich}. In fact, more can be said on the stability behavior of~\eqref{squareduffing} in a neighborhood of $(\delta,\gamma)=(0,0)$.
Even if we initially assumed that $\gamma\geq0$, we may give a full description of what happens when $\gamma<0$.

\begin{theorem}\label{striscia}
	The trivial solution of~\eqref{squareduffing} is:
	$$
	\mbox{stable if}\quad-\frac{\delta^2}{2}<\gamma<1\quad\forall\delta\in\R\, ,\qquad
	\mbox{unstable if}\quad\gamma<-\frac{\delta^2}{2}\quad\forall\delta\in\R\, .
	$$
\end{theorem}

\begin{figure}[ht]
	\begin{center}
		\includegraphics[height=0.3\textheight]{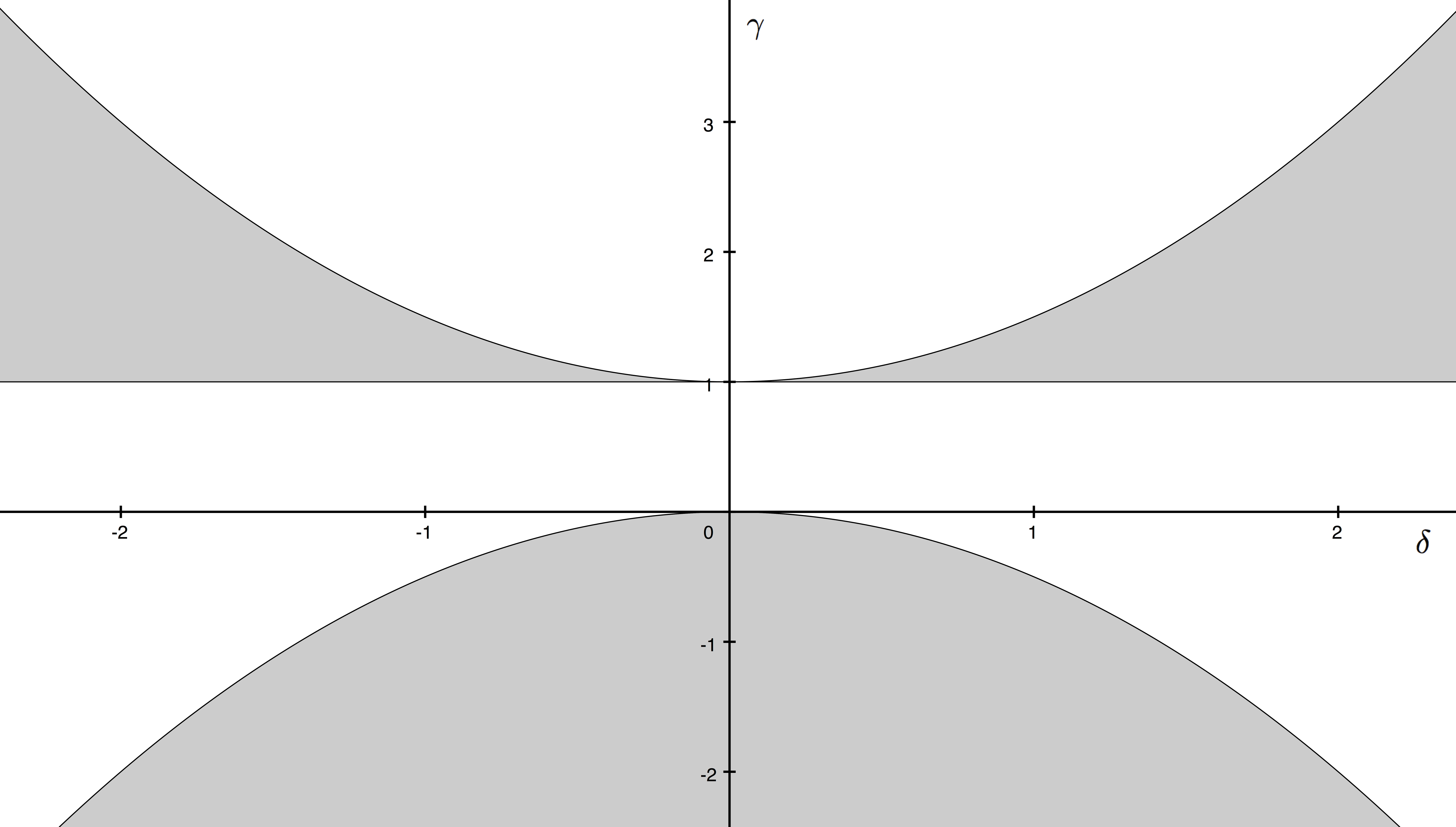}
		\caption{Stability regions (white) and resonance tongues (gray) for~\eqref{squareduffing}; see Theorems~\ref{exactsolutions}
and~\ref{striscia}.}\label{zonevere}
	\end{center}
\end{figure}

Also the proof of this result is given in Section~\ref{proofs}: in particular, we use there the fact that, if $\gamma\le-\delta^2$, then $\gamma+y(t)^2\le0$
and the instability is trivial. A picture of the stability region described by Theorem~\ref{striscia} is shown in Figure~\ref{zonevere}, where we have
included negative values for both $\gamma$ and $\delta$.

Unfortunately, we could not find an explicit representation of the remaining resonant lines, namely the boundaries of the resonance tongues
$U_\ell$ for $\ell\ge2$. Thus, we now use the Burdina criterion to determine sufficient conditions for the stability of the trivial solution of
\eqref{squareduffing}. To this end, we introduce the function
$$\Phi(\delta,\gamma):=2\sqrt{2}\int_0^{\pi/2}\sqrt{\tfrac{\gamma+\delta^2\sin^2\theta}{2+\delta^2+\delta^2\sin^2\theta}}\, d\theta
\qquad\forall\delta,\gamma>0\, .$$
Note that $\Phi$ is positive and smooth; this function may also be expressed in terms of elliptic integrals.\par
In Section~\ref{proofs} we will prove the following statement.
\begin{theorem}\label{primo}
	Let $y$ be the solution of the Duffing equation \eqref{ODE} with initial conditions \eqref{alphabeta}. If there exists $\ell\in\N$ such that
	\neweq{primavera}
	\log\left(1+\frac{\delta^2}{\gamma}\right)<2\cdot\min\Big\{\Phi(\delta,\gamma)-\ell\pi\, ,\, (\ell+1)\pi-\Phi(\delta,\gamma)\Big\}\, ,
	\endeq
	then the trivial solution of \eqref{squareduffing} is stable.
\end{theorem}

The trick in the proof is to use the Burdina criterion {\em without} using the explicit form of the solution of
\eq{ODE}. In particular, Theorem \ref{primo} has the following elegant consequence:
\begin{corollary}\label{corparabola}
	Let $y$ be the solution of the Duffing equation \eqref{ODE} with initial conditions \eqref{alphabeta}. Then
	the trivial solution of
	$$\ddot{\xi}(t)+\Big(2+\delta^2+y(t)^2\Big)\xi(t)=0$$
	is stable.
\end{corollary}

The proof follows directly from Theorem~\ref{primo} (case $\ell=1$) by observing that $\Phi(\delta,2+\delta^2)=\sqrt{2}\pi$ and
$$
\log\left(1+\frac{\delta^2}{2+\delta^2}\right)<\log 2<2\pi(\sqrt2 -1)=2\pi\cdot\min\Big\{\sqrt2 -1\, ,\, 2-\sqrt2 \Big\}\qquad\forall\delta>0\, .
$$

In order to obtain a more precise picture of the resonant tongues other than the first, we use Theorem \ref{primo} to deduce some information on their behavior as $\delta\to0$.
\begin{theorem}\label{asymp}
Let $\ell\in\N$ $(\ell\ge2)$ and let $U_\ell$ be the resonant tongue emanating from $(\delta,\gamma)=(0,\ell^2)$. If $(\delta,\gamma)\in U_\ell$, then
	\neweq{twoparabolas}
	\ell^2+\left(\frac{3\ell^2}{4}-\frac12 -\frac{1}{\pi\ell}\right)\delta^2+O(\delta^4)\le\gamma\le
	\ell^2+\left(\frac{3\ell^2}{4}-\frac12 +\frac{1}{\pi\ell}\right)\delta^2+O(\delta^4)\quad\mbox{as }\delta\to0\, .
	\endeq
	Therefore, for all $\gamma\ge0$ (with $\gamma\neq1$) there exists $\delta_\gamma>0$ such that if $0<\delta<\delta_\gamma$ and if $y$ is the solution of
	\eqref{ODE}-\eqref{alphabeta}, then the trivial solution of \eqref{squareduffing} is stable.
\end{theorem}

\begin{figure}[ht]
	\begin{center}
		\includegraphics[height=0.22\textheight]{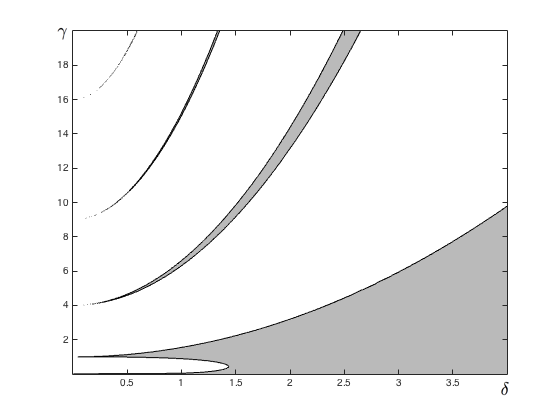}\qquad\qquad\includegraphics[height=0.22\textheight]{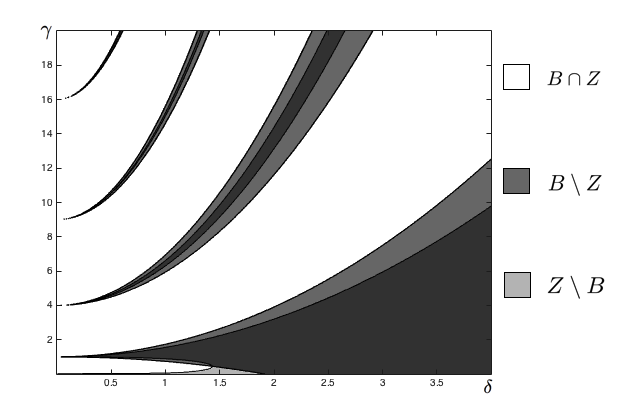}
		\caption{Stability regions (white) for \eq{squareduffing} obtained numerically through the
			Burdina criterion (left) and comparison with the Zhukovskii criterion (right): in the legend, the set $Z$ contains
the points $(\delta,\gamma)$ which satisfy the Zhukovskii criterion, and $B$ contains the points which satisfy the Burdina criterion.}\label{num_B}
	\end{center}
\end{figure}

A further remark about the comparison between the three criteria reported in Proposition \ref{lyapzhu} is in order. In the left picture of Figure \ref{num_B}
we plot the (white) stability regions obtained numerically through the Burdina criterion, see \eq{primavera}. The second lowest stable (white) region contains
the parabola $\gamma=2+\delta^2$. In the right picture of Figure \ref{num_B} we see that the Burdina criterion performs better than the Zhukovskii criterion,
at least in the region where $\gamma>1$: except for the first stability region, the regions which result stable for the Zhukovskii criterion are strictly included in the ones determined by the Burdina criterion. On the other hand, back to the left picture, we see that the Burdina criterion
does not ensure stability in a neighborhood of $(0,0)$; in this case the Zhukovskii criterion performs slightly better. But the $L^2$ generalized Lyapunov
criterion by Li-Zhang is the one having the better performance in the region $\gamma<1$: a hint of this fact is highlighted in the proof of
Theorem \ref{intornorigine}.

\begin{figure}[ht]
	\begin{center}
		\includegraphics[height=0.4\textheight]{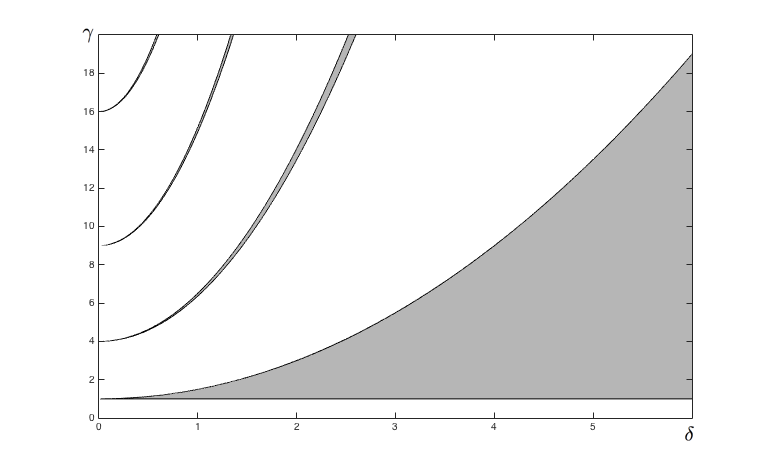}
		\caption{Stability regions (white) for \eq{squareduffing} obtained numerically with the use of the monodromy matrix.}\label{realplot}
	\end{center}
\end{figure}

Theorem \ref{primo} only gives sufficient conditions for the stability and the stability regions are considerably wider than the those defined
by \eq{primavera}. In order to have a more precise picture of the resonant tongues, we proceed numerically. In Figure \ref{realplot} we plot
the obtained stability (white) and instability (gray) regions for \eq{squareduffing}. It turns out that the resonant tongues $U_\ell$ (with $\ell\ge2$)
are extremely narrow, see \eq{twoparabolas}.
The boundaries of $U_\ell$ (resonant lines) are found numerically by computing the trace of the monodromy matrix of the Hill
equation \eqref{squareduffing} and by plotting its level lines $\pm2$.
To do so, we used the software Matlab. The first step was to introduce a grid of values of $\delta$ and $\gamma$: for each value of $\delta$, with the built-in function \texttt{ellipj}, which gives the elliptic Jacobi functions, we defined a Matlab function which computed the solution $y_\delta$
of the Duffing equation with initial conditions \eq{alphabeta}, using~\eqref{soluzione_duffing}. Then, by using the Matlab ODE solver \texttt{ode45}, we integrated the principal fundamental matrix at $t=0$ in the interval $[0,T(\delta)]$ for every point $(\delta,\gamma)$. This method introduced some errors,
both in the computation of the Jacobi cosine elliptic function and in the integration of the ODE. For this reason, and because the resonant tongues are very
thin in a neighborhood of $\delta=0$, the plot was quite difficult. In order to give a more clear representation of the tongues, we plotted the level
lines $\pm1.98$ instead of $\pm2$ of the trace of the monodromy matrix. The result is shown in Figure \ref{realplot}.

\section{Instabilities in a nonlinear nonlocal beam equation}\label{appl}

\subsection{Nonlinear modes}

In this section we define the nonlinear modes of \eq{truebeam} and what we mean by stability.
We start by seeking particular solutions of \eq{truebeam} by separating variables, that is, in the form
\neweq{form}
u_m(x,t)=\Theta_m(t)\sin(mx)\qquad(m\in\N)\, .
\endeq
A simple computation shows that the function $\Theta_m$ satisfies
\neweq{ODED}
\ddot{\Theta}_m(t)+m^4\Theta_m(t)+m^4\Theta_m(t)^3=0\qquad(t>0)\, ,
\endeq
which is a time-scaled version ($t\mapsto m^2t$) of the Duffing equation \eq{ODE}.\par
We call a function $u_m$ in the form \eq{form} a $m$-th {\bf nonlinear mode} of \eq{truebeam}. Note that for any $m$ there exist infinitely many
$m$-th nonlinear modes depending on the initial value, they are not proportional to each other and they have different periodicity.
Their shape is described by the solutions $\Theta_m$ of \eq{ODED}: for this reason, with an abuse of language, we also call $\Theta_m$ a nonlinear mode
of \eq{truebeam}.\par
We are interested in studying the stability of the nonlinear modes. To this end, we consider solutions of \eq{truebeam} in the form
\neweq{form2}
u(x,t)=w(t)\sin(mx)+z(t)\sin(nx)
\endeq
for some integers $n,m\ge1$, $n\neq m$. After inserting \eq{form2} into \eq{truebeam} we reach the following (nonlinear) system of ODE's:
\neweq{cw}
\left\{\begin{array}{l}
\ddot{w}(t)+m^4w(t)+m^2\big(m^2w(t)^2+n^2z(t)^2\big)w(t)=0\, ,\\
\ddot{z}(t)+n^4z(t)+n^2\big(m^2w(t)^2+n^2z(t)^2\big)z(t)=0\ ,
\end{array}\right.
\endeq
to which we associate some initial conditions
\neweq{initialsyst}
w(0)=w_0\, ,\ \dot{w}(0)=w_1\, ,\quad z(0)=z_0\, ,\ \dot{z}(0)=z_1\, .
\endeq
The constant energy of the Hamiltonian system \eq{cw} is given by
\begin{eqnarray}
\E(w_0,w_1,z_0,z_1)&=&\frac{\dot{w}^2}{2}+\frac{\dot{z}^2}{2}+m^4\frac{w^2}{2}+n^4\frac{z^2}{2}+
\frac{(m^2w^2+n^2z^2)^2}{4}\notag \\
&\equiv&\frac{w_1^2}{2}+\frac{z_1^2}{2}+m^4\frac{w_0^2}{2}+n^4\frac{z_0^2}{2}+\frac{(m^2w_0^2+n^2z_0^2)^2}{4}\, .\label{uffa}
\end{eqnarray}

Note that if $z_0=z_1=0$, then the solution of \eq{cw}-\eq{initialsyst} is $(w,z)=(\Theta_m,0)$ where $\Theta_m$ is a nonlinear mode, namely a
solution of the Duffing equation \eq{ODED}. In order to analyze the stability of the nonlinear mode $\Theta_m$ with respect to the nonlinear mode
$\Theta_n$ we argue as follows. We take initial data in \eq{initialsyst} such that
$$
0<|z_0|+|z_1|\ll|w_0|+|w_1|\, .
$$
This means that the energy \eq{uffa} is initially almost totally concentrated on the $m$-th mode, that is, $\E(w_0,w_1,z_0,z_1)\approx\E(w_0,w_1,0,0)$.
We then wonder whether this remains true for all time $t>0$ for the solution of \eq{cw}. To this end, we introduce the following notion of stability.

\begin{definition}\label{defstabb}
The mode $\Theta_m$ is said to be {\bf linearly stable} ({\bf unstable}) with respect to the $n$-th mode $\Theta_n$ if $\xi\equiv0$ is a stable
(unstable) solution of the linear Hill equation
\neweq{hill2}
\ddot{\xi}(t)+\Big(n^4+m^2n^2\Theta_m(t)^2\Big)\xi(t)=0\, ,\quad \forall t\, .
\endeq
\end{definition}

There exist also stronger definitions of stability which, in some cases, can be shown to be equivalent; see e.g.\ \cite{ghg}. For nonlinear PDE's
such as \eq{truebeam}, Definition \ref{defstabb} is sufficiently precise to characterize the instabilities of the nonlinear modes
of the equation, see \cite{bbfg,bfg1,bergaz,fergazmor} and also the example at the end of the next subsection.\par
The relevant parameter for the stability of the trivial solution of \eq{hill2} turns out to be
\[\omega:=\frac{n^2}{m^2}\, .\]
Note that, if $\Theta_m(t)$ is a solution of~\eqref{ODED} with initial conditions
\neweq{initialOmega}
\Theta_m(0)=\delta\quad\mbox{and}\quad\dot{\Theta}_m(0)=0\, ,
\endeq
then $\Theta_\omega(t)=\Theta_m(t/mn)$ solves the following time-scaled version of the Duffing equation \eq{ODE}:
\begin{equation}\label{duffomega}
\ddot{\Theta}_\omega(t)+\frac{1}{\omega}\Theta_\omega(t)+\frac{1}{\omega}\Theta_\omega(t)^3=0,
\end{equation}
with the same initial conditions. From Burgreen~\cite{burg}, we know that
\begin{equation}\label{solomega}
\Theta_\omega(t)=\delta\cn\Bigg[t\,\sqrt{\frac{1+\delta^2}{\omega}},\frac{\delta}{\sqrt{2(1+\delta^2)}}\Bigg].
\end{equation}
Moreover, its period is given by
\begin{equation*}\label{Tomega}
T_\omega(\delta)=4\,\sqrt{\frac{\omega}{1+\delta^2}}\,K\bigg(\frac{\delta}{\sqrt{2(1+\delta^2)}}\bigg),
\end{equation*}
so that $\lim_{\delta\to 0}T_\omega(\delta)=2\pi\sqrt{\omega}$. The counterpart of~\eqref{ancoraenergia} reads
\begin{equation}
2\dot{\Theta}_\omega^2=\frac{1}{\omega}(2+\delta^2+\Theta_\omega^2)(\delta^2-\Theta_\omega^2)\qquad\forall t\, .
\end{equation}
Through the time scaling $t\mapsto\frac{t}{mn}$,~\eqref{hill2} is equivalent to the Hill equation
$$
\ddot{\xi} + \bigg(\frac{n^2}{m^2} + \Theta_m\bigg(\frac{t}{mn}\bigg)^2\bigg)\xi=0\, ,
$$
that is,
\begin{equation}\label{hill_scalata}
\ddot{\xi}(t)+(\omega+\Theta_\omega(t)^2)\xi(t)=0\, ,
\end{equation}
where $\Theta_\omega$ is given in \eq{solomega}.

\subsection{Stability of the nonlinear modes}

We study here in detail the stability of the nonlinear modes of \eq{truebeam}, according to Definition \ref{defstabb}.
In particular, we prove the following statement, left open in \cite{bbfg}.

\begin{theorem}\label{facile}
If $m>n$ then the mode $\Theta_m$ is linearly stable with respect to the mode $\Theta_n$, independently of the initial value $\delta\neq0$ in \eqref{initialOmega}.
\end{theorem}

The proof of this result follows from a more general statement, see Theorem \ref{primalingua} below, therefore we omit it.\par
The situation is by far more complicated when $m<n$, several subcases have to be distinguished.
From~\cite[Theorem 15]{bbfg} and~\cite[Theorem 1.1]{cazw2}, we learn that an asymptotic representation of the resonant tongues as $\delta\to\infty$
in \eq{initialOmega} can be given. Let us define the two sets
\begin{equation}\label{strisceinfinito}
I_U:=\bigcup_{k\in\N}\Big((k+1)(2k+1),(k+1)(2k+3)\Big)\qquad
I_S:=\bigcup_{k\in\N}\Big(k(2k+1),(k+1)(2k+1)\Big).
\end{equation}

Note that $\overline{I_S\cup I_U}=[0,\infty)$. Then, by using the very same method as in \cite{cazw2}, the following result was obtained in \cite{bbfg}:

\begin{proposition}\label{known}
Let $I_S$ and $I_U$ be as in \eqref{strisceinfinito} and let $\Theta_\omega$ be as in \eqref{solomega}. For every $\omega>0$,
there exist $\bar{\delta}_\omega>0$ such that, for all $\delta>\bar{\delta}_\omega$:\par\noindent
$(i)$ if $\omega\in I_U$, then the trivial solution of equation~\eqref{hill_scalata} is unstable and $\Theta_m$ is linearly unstable
with respect to $\Theta_n$;\par\noindent
$(ii)$ if $\omega\in I_S$, then the trivial solution of equation~\eqref{hill_scalata} is stable and $\Theta_m$ is linearly stable
with respect to $\Theta_n$.
\end{proposition}

Similarly to what we did in Section~\ref{duffsqMOD}, we apply the Burdina criterion to equation~\eqref{hill_scalata}. To this end, we define the function
\begin{equation}\label{Psi}
\Psi(\delta,\omega)=2\sqrt{2\,\omega}\int_0^{\pi/2}\sqrt{\tfrac{\omega+\delta^2\sin^2\theta}{2+\delta^2+\delta^2\sin^2\theta}}\, d\theta
\qquad\forall\delta,\omega>0\, .
\end{equation}
Note that the map $\delta\mapsto\Psi(\delta,\omega)$ is strictly decreasing for all $\omega\ge1$ and that
$$\Psi(0,\omega)=\pi\omega\, ,\qquad\lim_{\delta\to\infty}\Psi(\delta,\omega)=
2\sqrt{2\,\omega}\int_0^{\pi/2}\frac{\sin\theta}{\sqrt{1+\sin^2\theta}}\, d\theta=\pi\, \sqrt{\frac{\omega}{2}}\, .$$

By applying the Burdina criterion, see Proposition \ref{lyapzhu} $(iii)$, we obtain the following sufficient condition for the stability
of the trivial solution of~\eqref{hill_scalata}.

\begin{theorem}\label{burdinaomega}
Let $\Theta_\omega$ be as in \eqref{solomega} and let $\Psi$ be as in \eqref{Psi}. If there exists $\ell\in\N$ such that
\begin{equation*}\label{fantasia}
\log\Big(1+\frac{\delta^2}{\omega}\Big)<2\cdot\min\Big\{\Psi(\delta,\omega)-\ell\pi\, ,\, (\ell+1)\pi-\Psi(\delta,\omega)\Big\},
\end{equation*}
then the trivial solution of~\eqref{hill_scalata} is stable and $\Theta_m$ is linearly stable with respect to $\Theta_n$.
\end{theorem}

\begin{figure}[ht]
	\begin{center}
		\includegraphics[height=0.4\textheight]{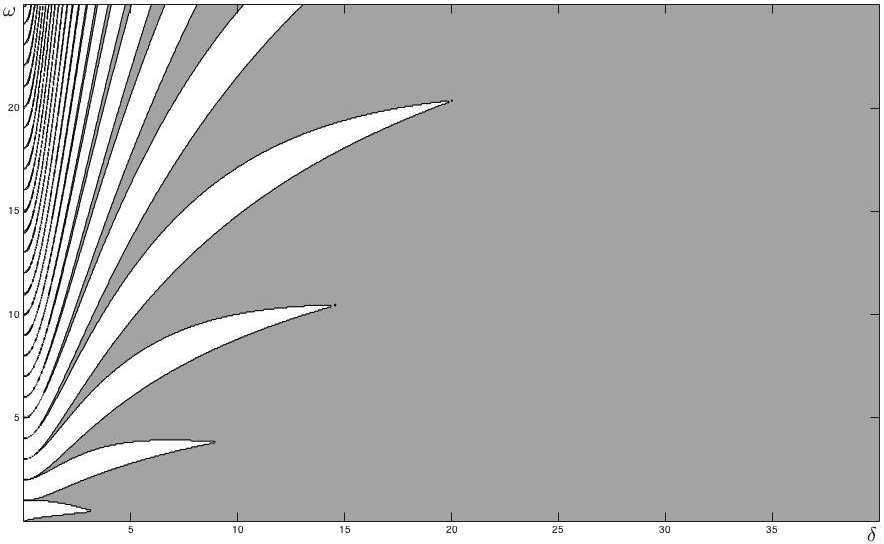}
		\caption{Stability regions (white) obtained with the sufficient condition of Theorem \ref{burdinaomega}.}\label{burdina2}
	\end{center}
\end{figure}

In Figure \ref{burdina2} we display the regions described by Theorem \ref{burdinaomega}.

In particular, Theorem~\ref{burdinaomega} enables us to determine the behavior of the resonant tongues as $\delta\to0$.

\begin{corollary}\label{asympomega}
	Let $\ell\in\N$ $(\ell\ge2)$ and let $U_\ell$ be the resonant tongue of \eqref{hill_scalata} emanating from $(\delta,\omega)=(0,\ell)$.
    If $(\delta,\omega)\in U_\ell$, then
	$$
	\ell+\left(\frac{3\ell}{8}-\frac{1}{4} -\frac{1}{2\pi\ell}\right)\delta^2+O(\delta^4)\le\omega\le
	\ell+\left(\frac{3\ell}{8}-\frac{1}{4} +\frac{1}{2\pi\ell}\right)\delta^2+O(\delta^4)\quad\mbox{as }\delta\to0\, .
	$$
\end{corollary}

The next result, whose proof is given in Section~\ref{proofs}, characterizes the first resonance tongue.

\begin{theorem}\label{primalingua}
The resonant tongue $U_1$ of \eqref{hill_scalata} emanating from $(\delta,\omega)=(0,1)$ is given by
$$
U_1 =\big\{(\delta,\omega)\in\mathbb{R}_+^2:\, 1<\omega<\phi(\delta)\big\}
$$
where $\phi:\mathbb{R}_+\to\mathbb{R}_+$ is a continuous function such that
$$
\phi(\delta)>1\ \forall\delta>0\, ,\qquad\phi(\delta)\le1+\bigg(\frac{1}{8}+\frac{1}{2\pi}\bigg)\delta^2+O(\delta^4)\mbox{ as }\delta\to0\, ,\qquad
\lim_{\delta\to\infty}\phi(\delta)=3\, .
$$
Moreover, the strip $(0,\infty)\times(0,1)$ is a stability region of the $(\delta,\omega)$-plane.
\end{theorem}

A straightforward consequence of Corollary \ref{asympomega} and Theorem \ref{primalingua} reads

\begin{corollary}\label{newcor}
For all $\omega>0$ (with $\omega\neq1$) there exists $\delta_\omega>0$ such that if $0<\delta<\delta_\omega$ and if $\Theta_\omega$ is as
in \eqref{solomega}, then the trivial solution of~\eqref{hill_scalata} is stable. Therefore, for any integers $m\neq n$, there exists
$\widehat{\delta}=\widehat{\delta}(n/m)>0$ such that the nonlinear mode $\Theta_m$, solution of \eqref{ODED}-\eqref{initialOmega},
is linearly stable with respect to $\Theta_n$ whenever $0<|\delta|<\widehat{\delta}$; moreover, $\widehat{\delta}=+\infty$ if and only if $n<m$.
\end{corollary}

\begin{figure}[ht]
	\begin{center}
		\includegraphics[height=0.4\textheight]{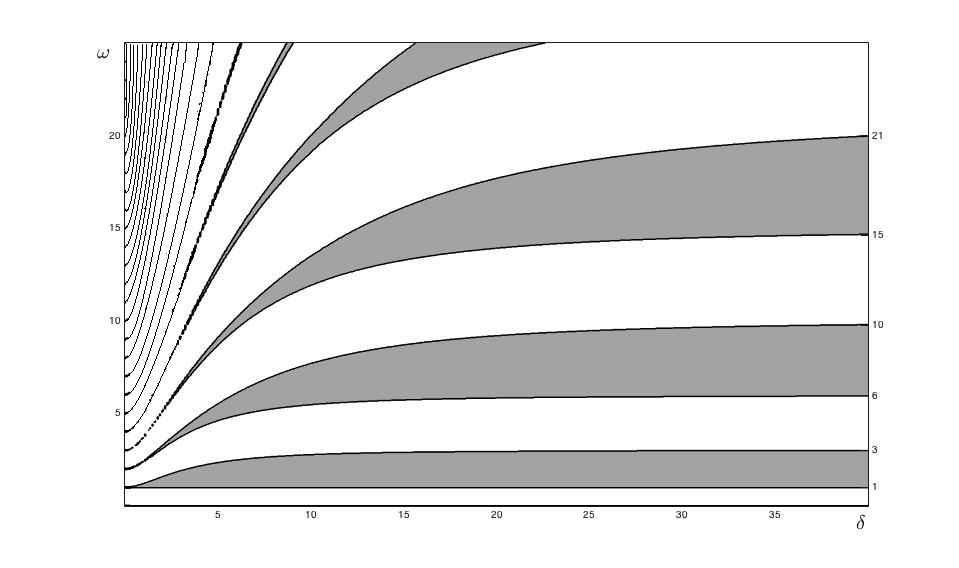}
		\caption{Stability regions (white) for~\eq{hill_scalata} obtained numerically with the use of the monodromy matrix.}\label{monodromymnd1}
	\end{center}
\end{figure}

By proceeding as in Section~\ref{duffsqMOD}, we numerically obtained a complete picture of the resonance tongues. We fixed a value of $m$
in~\eqref{hill2} and we integrated system~\eqref{ODED}-\eqref{hill2} allowing $n$ to take non integer values, so that also the squared ratio
$\omega$ in~\eqref{hill_scalata} could take any real positive value. The resonant lines were obtained by computing the trace of the
monodromy matrix of the Hill equation~\eqref{hill2} for different values of $\delta$ and $\omega$: since the resonance tongues are very narrow for small
$\delta$, the plot of the level curves $\pm2$ in the $(\delta,\omega)$-plane was quite difficult and we obtained an approximate picture
of the stability regions by plotting instead the level curves $\pm1.98$. The result is shown in Figure~\ref{monodromymnd1}.

Due to the asymptotic behavior of the stability regions stated in Proposition~\ref{known}, the resonance tongues are much wider for large
values of $\delta$ than in a neighborhood of the $\omega$ axis. For this reason, if we choose a value $\omega$ of the squared ratio $n^2/m^2$
and we increase $\delta$ starting from the point $(0,\omega)$, we cross some tiny resonance tongues before reaching the final
stability or instability region, as described in Proposition~\ref{known}.
For a given squared ratio $\omega=n^2/m^2$ of spatial frequencies, in Table \ref{crossings} we quote the (minimum) number of crossing of the resonance
lines, starting from $\delta=0$ up to $\delta\to\infty$. Since by Corollary \ref{newcor} we know that for $\delta\to0^+$ we always start in a stability
region, if this number is even (resp. odd) then we end up in a stability (resp. instability) region when $\delta\to\infty$.
\begin{table}[ht]
\begin{center}
\begin{tabular}{|c|c|c|c|c|c|c|c|c|c|}
\hline
$\omega\in$ & $(0,1)$ & $(1,2]$ & $(2,3)$ & $\{3\}$ & $(3,4]$ & $(4,5]$ & $(5,6)$ & $\{6\}$ & $(6,7]$\\
\hline
crossing & $0$ & $1$ & $3$ & $2$ & $4$ & $6$  & $8$ & $7$ & $9$\\
\hline
\end{tabular}
\caption{Number of resonance lines crossed before reaching the final stability/instability region.}\label{crossings}
\end{center}
\end{table}

Let us emphasize that different couples $(m,n)$ may have the very same stability behavior. Fix some $(m,n)$ and consider all the couples
$(km,kn)$ for any integer $k$. Then the corresponding $\omega$ is the same and hence, also the stability behavior. What changes is the time scaling:
by applying the scaling $t\mapsto kt$, we recover the very same system \eq{cw} and the occurrence of a possible instability (as in the plots of
Figure \ref{plots3} below) will appear delayed in time. Whence, what really counts for the stability of nonlinear modes is the {\em ratio}
$\omega=n^2/m^2$ of spatial frequencies.\par
Let us conclude with a particular example which shows how the instability for \eq{cw}-\eq{initialsyst} appears.
We take $m=1$ and $n=2$, so that $\omega=4$ and we have the following consequence of Theorem \ref{burdinaomega}.

\begin{corollary}\label{cor4}
If one of the following two facts holds
$$\log\left(1+\frac{\delta^2}{4}\right)<2\cdot\min\big\{\Psi(\delta,4)-2\pi,3\pi-\Psi(\delta,4)\big\}\, ,$$
$$\log\left(1+\frac{\delta^2}{4}\right)<2\cdot\min\big\{\Psi(\delta,4)-3\pi,4\pi-\Psi(\delta,4)\big\}\, ,$$
then the mode $\Theta_1$ is linearly stable with respect to the mode $\Theta_2$.
\end{corollary}

Numerically, one sees that Corollary \ref{cor4} guarantees the linear stability of the mode $\Theta_1$ with respect to the mode $\Theta_2$ whenever
\neweq{interval}
\delta\in(0,1.167)\cup(1.277,2.63)\, ,
\endeq
while Proposition \ref{known} guarantees linear stability for $\delta>\overline{\delta}$ for a sufficiently large $\overline{\delta}$.
Therefore, the instability range for the initial semi-amplitude $\delta$ lies in the complement of \eq{interval}. According to Figure
\ref{monodromymnd1} we expect it to be composed by two disconnected intervals. In order to find it with precision, we numerically plotted
the solution of \eq{hill_scalata} with $\Theta_\omega$ defined in \eq{solomega} and $\omega=4$. We could observe instability, that is
exponential-like blow up of the solution towards $\pm\infty$ when $t\to+\infty$, for $\delta\in(2.93,3.45)$ which lies in the complement of
\eq{interval} and belongs to the resonance tongue emanating from $(0,2)$ in the $(\delta,\omega)$-plane. This also means that $\overline{\delta}\approx3.45$.
Then we increased $\delta$ with step $10^{-4}$
but we detected no instability for $\delta\in(1.167,1.277)$ which contains the intersection of the resonant tongue emanating
from $(\delta,\omega)=(0,3)$ with the line $\omega=4$: this means that this region is extremely thin and/or that
the modulus of the trace of the monodromy matrix only slightly exceeds 2. Overall, this means that these instabilities are irrelevant:
on the one hand, they have very small probability to appear, on the other hand, even if they do appear they only give rise to tiny transfers
of energy. For any $\omega$ one then expects to view ``true'' instabilities only for large values of $\delta$.\par
We also plotted the solution of \eq{cw}-\eq{initialsyst}: we considered the system
\neweq{cw12}
\left\{\begin{array}{l}
\ddot{w}(t)+w(t)+\big(w(t)^2+4z(t)^2\big)w(t)=0\, ,\\
\ddot{z}(t)+16z(t)+4\big(w(t)^2+4z(t)^2\big)z(t)=0\, ,\\
w(0)=\delta\, ,\ z(0)=10^{-3}\delta\, ,\ \dot{w}(0)=\dot{z}(0)=0\, .
\end{array}\right.
\endeq
In Figure \ref{plots3} we display some of the pictures that we obtained.
\begin{figure}[ht]
\begin{center}
\includegraphics[height=44mm,width=84mm]{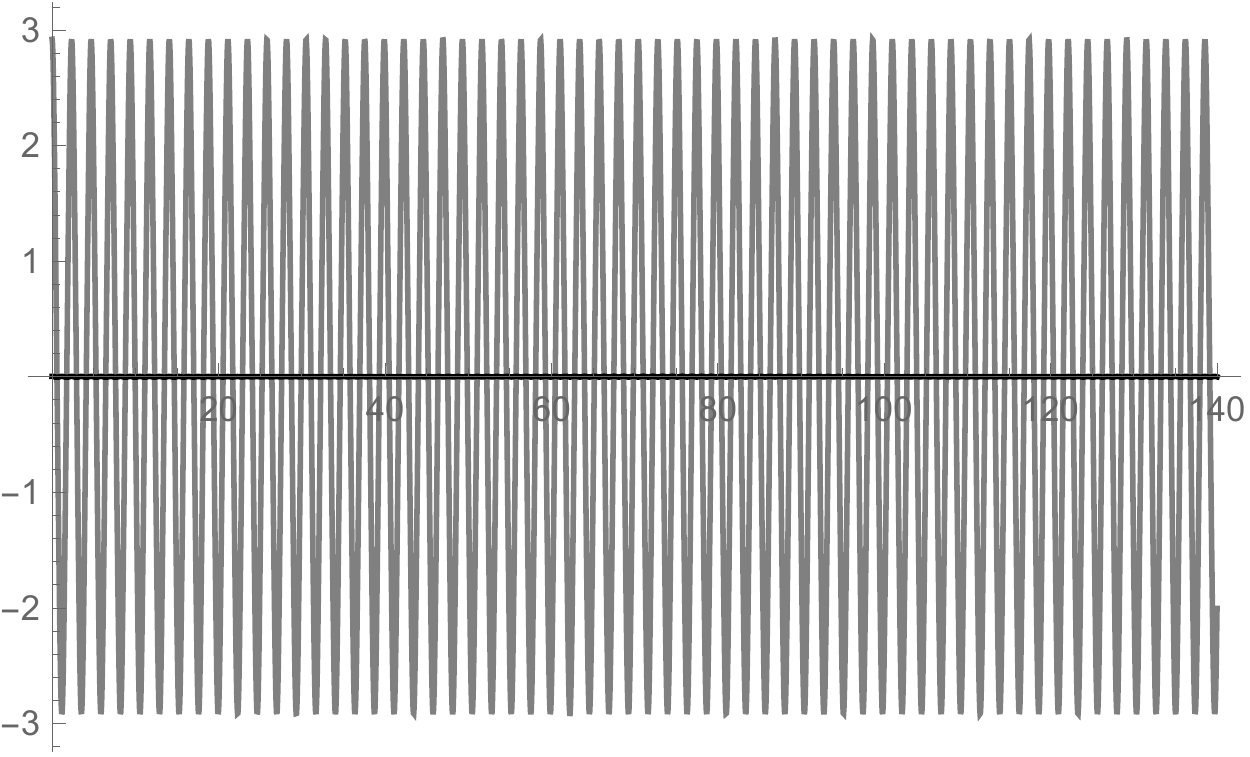}\quad\includegraphics[height=44mm,width=84mm]{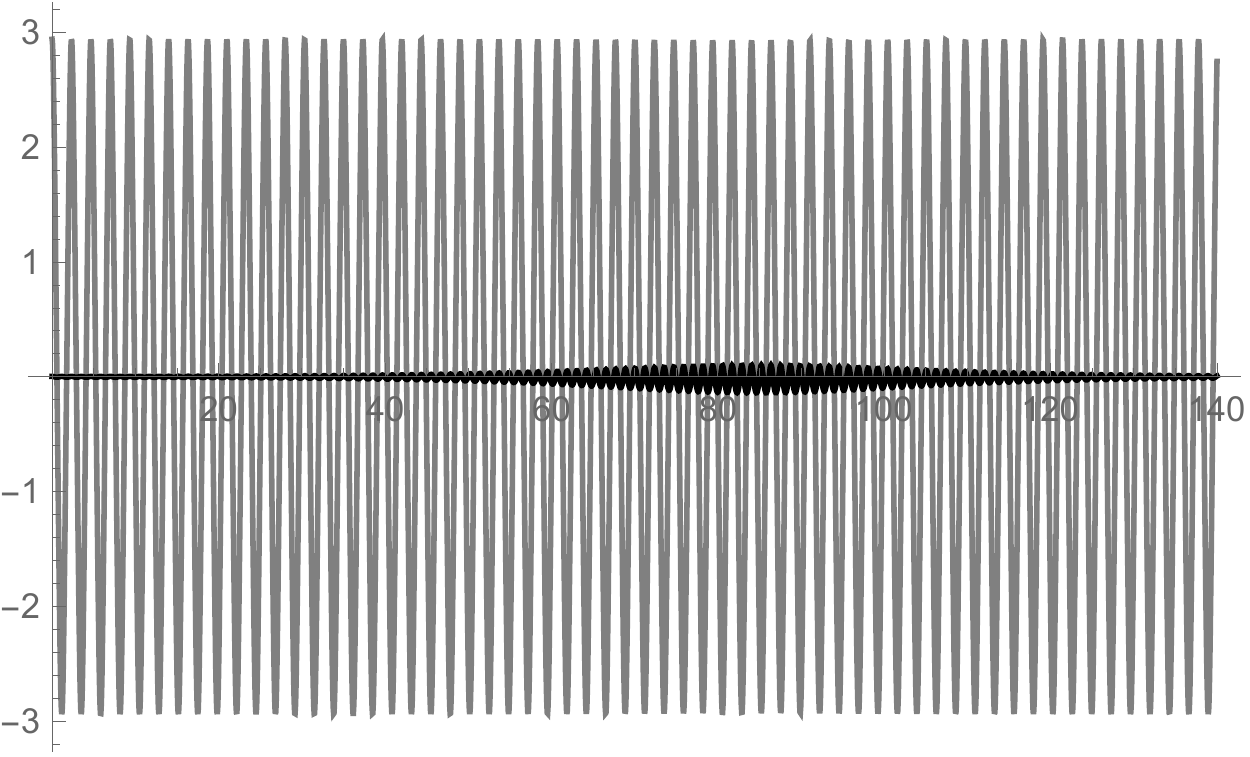}

\includegraphics[height=44mm,width=84mm]{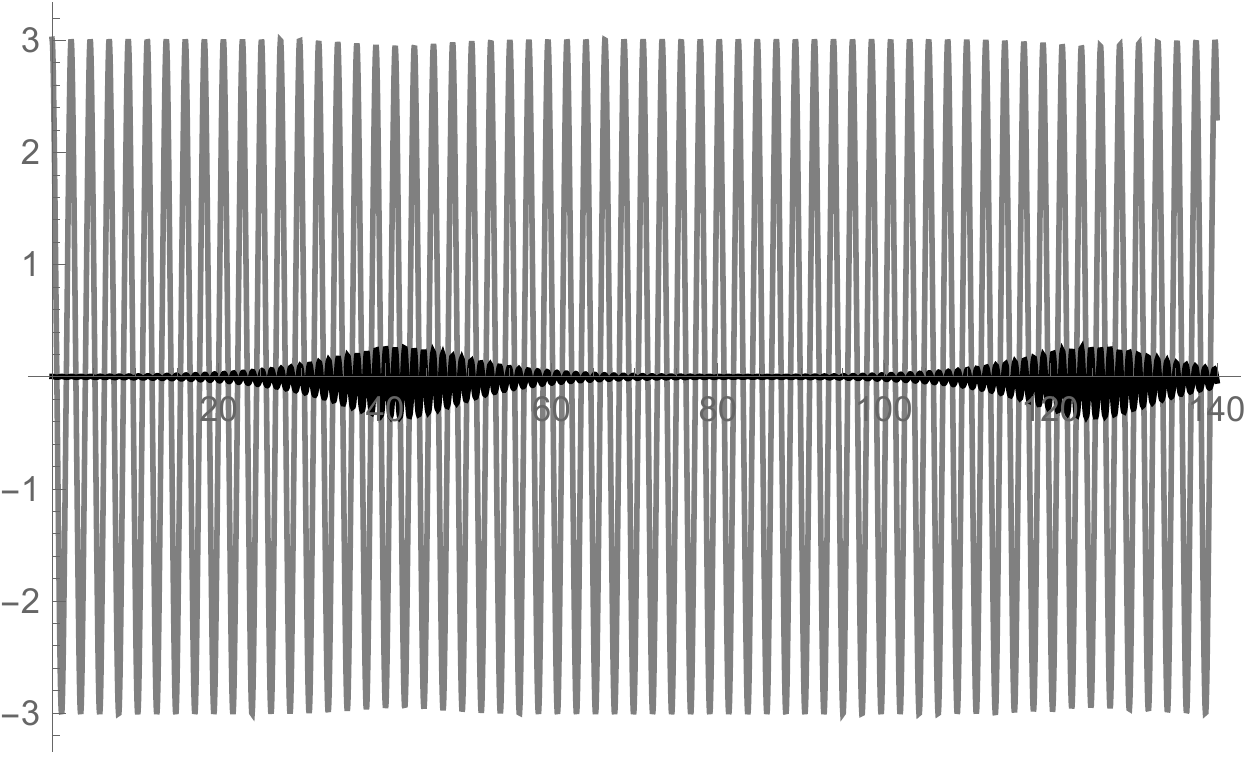}\quad\includegraphics[height=44mm,width=84mm]{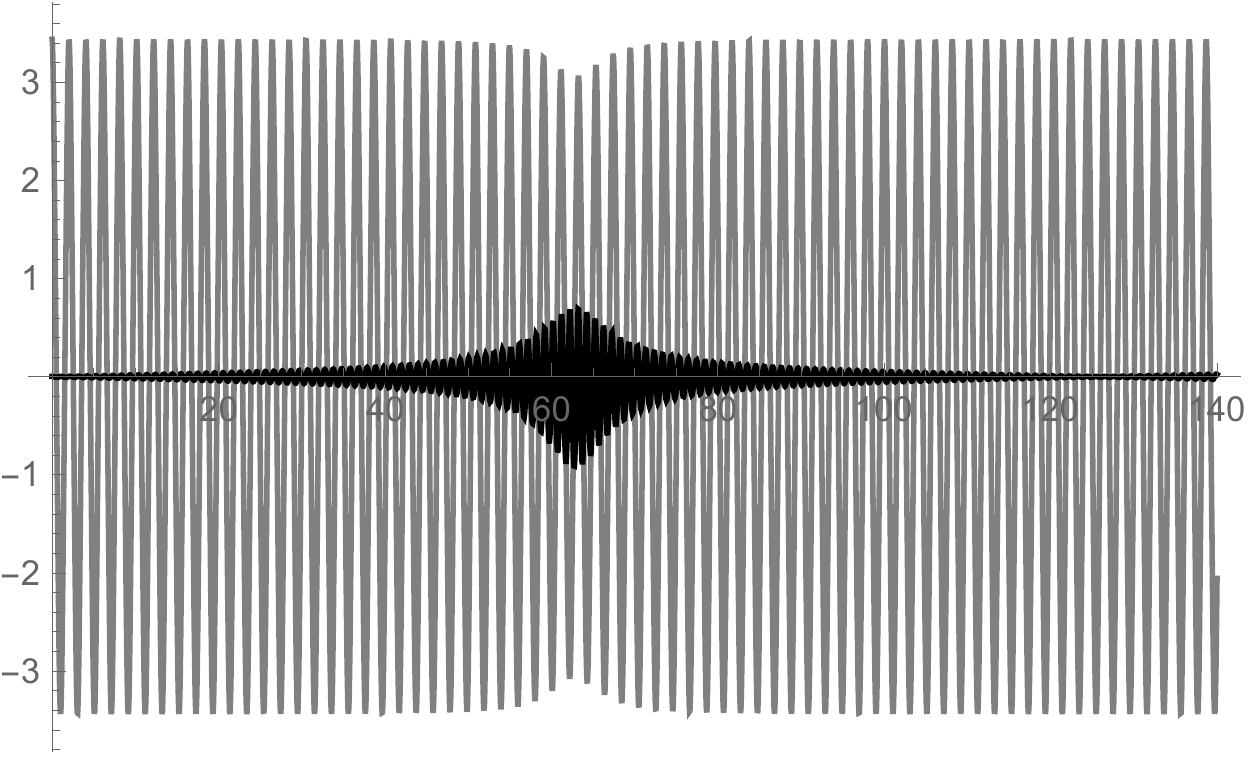}
\caption{The solutions $w$ (gray) and $z$ (black) of \eq{cw12} for $\delta=2.92$, $\delta=2.94$, $\delta=3.01$, $\delta=3.44$.}\label{plots3}
\end{center}
\end{figure}
It appears clearly that, in the first picture, $z$ remains small and no instability occurs. In the second picture, $z$ suddenly grows up, which means
instability. In the third picture the same phenomenon is accentuated while in the fourth picture it is delayed, thereby getting ready to return in a stability regime
for $\delta>3.454$. This is the pattern that can be observed in \eq{cw}-\eq{initialsyst}, for any
value of $\omega=n^2/m^2$. This perfect coincidence with the linear stability justifies Definition \ref{defstabb} and the use of linearization.

\section{Proofs}\label{proofs}

Since we need a formula stated therein, we prove Theorem~\ref{primo} first.

\noindent{\bf Proof of Theorem~\ref{primo}.} We may assume that $y(0)=-\delta<0$ so that $\dot y>0$ in $(0,T/2)$. Then we get
\begin{eqnarray}
\int_0^{T/2}\sqrt{\gamma+y(t)^2}\, dt &=& 2\int_0^{T/4}\sqrt{\gamma+y(t)^2}\, dt \notag \\
\mbox{by \eq{derivata} \ } &=& 2\sqrt{2}\int_0^{T/4}\sqrt{\tfrac{\gamma+y(t)^2}{(2+\delta^2+y(t)^2)(\delta^2-y(t)^2)}}\, \dot{y}(t)\, dt \notag \\
\mbox{using $y(t)=\delta\sin\theta$ \ } &=& \Phi(\delta,\gamma)\, .\label{super}
\end{eqnarray}
Moreover, if we put $p(t)=\gamma+y(t)^2$, we see that
$$\log\frac{\max p}{\min p}=\log\left(1+\frac{\delta^2}{\gamma}\right)\, .$$
Then, by Proposition \ref{lyapzhu} ($iii$), we infer that the trivial solution of \eqref{squareduffing} is stable whenever \eq{primavera} holds.\endproof

\noindent
{\bf Proof of Theorem~\ref{intornorigine}.} By using \eqref{TE} and by arguing as for~\eqref{super}, we obtain
\begin{align}
g(\delta)&:=\frac{T(\delta)^3}{8} \int_{0}^{T(\delta)/2} y(t)^4\,dt\notag\\
\ &=(2\sqrt{2})^4 \bigg(\int_{0}^{T(\delta)/4}\frac{y(t)^4\dot{y}(t)}{\sqrt{(2+\delta^2+y(t)^2)(\delta^2-y(t)^2)}}\,dt\bigg) \bigg(\int_{0}^{1}\frac{d\theta}{\sqrt{(2+\delta^2+\delta^2\theta^2)(1-\theta^2)}}\bigg)^3 \notag\\
\ &=64\bigg(\int_{0}^{\pi/2}\frac{\sin^4\alpha}{\sqrt{\frac{2}{\delta^2} + 1 + \sin^2\alpha}}\,d\alpha\bigg)\bigg(\int_0^{\pi/2}\frac{d\alpha}{\sqrt{\frac{2}{\delta^2} + 1 + \sin^2\alpha}}\bigg)^3.\notag
\end{align}
Hence, the map $\delta\mapsto g(\delta)$ is strictly increasing; furthermore, $g(0)=0$. Thus, for all $\delta>0$,
\begin{align}
0<g(\delta)<\lim_{\delta\to\infty}g(\delta)&=64\bigg(\int_{0}^{\pi/2}\frac{\sin^4\alpha}{\sqrt{1+\sin^2\alpha}}\,d\alpha\bigg)
\bigg(\int_{0}^{\pi/2}\frac{d\alpha}{\sqrt{1+\sin^2\alpha}}\bigg)^3	\notag\\
\ \mbox{using~\eqref{costantesigma}}\quad &=64\bigg(\int_{0}^{\pi/2}\frac{\sin^4\alpha}{\sqrt{1+\sin^2\alpha}}\,d\alpha\bigg)\sigma^3.	\label{bound_g}
\end{align}
Let us now put
$$
I := \int_{0}^{\pi/2}\frac{\sin^4\alpha}{\sqrt{1+\sin^2\alpha}}\,d\alpha
=\int_{0}^{\pi/2}\frac{\sin^2\alpha}{\sqrt{1+\sin^2\alpha}}\,d\alpha - \int_{0}^{\pi/2}\frac{\sin^2\alpha\cos^2\alpha}{\sqrt{1+\sin^2\alpha}}\,d\alpha\, ,
$$
so that, integrating by parts the second addend, we obtain
$$
I=\bigg(\int_{0}^{\pi/2}\frac{\sin^2\alpha}{\sqrt{1+\sin^2\alpha}}\,d\alpha\bigg)
+\bigg(\sigma-\int_{0}^{\pi/2}\frac{\sin^2\alpha}{\sqrt{1+\sin^2\alpha}}\,d\alpha -2I\bigg)\, ,
$$
which yields $I=\frac{\sigma}{3}$. Combined with \eqref{bound_g} and the monotonicity of $g$, this shows that $g(\delta)<\frac{64}{3}\sigma^4$ for
all $\delta>0$. We have so proved that
$$
\frac{T(\delta)^3}{8} \int_{0}^{T(\delta)/2} y(t)^4\,dt<\frac{64}{3}\sigma^4\qquad\forall\delta>0\, ,
$$
where $T(\delta)$ is the period of the solution $y(t)$ of the Duffing equation, as defined in Proposition~\ref{prop:sol_duffing}.
By recalling that the period of $y^2$ is $T(\delta)/2$, Proposition~\ref{lyapzhu} ($i$) applied to equation~\eq{squareduffing} and the latter inequality
ensure the stability of the trivial solution of~\eqref{squareduffing} if $\gamma=0$.\endproof

\noindent{\bf Proof of Theorem~\ref{exactsolutions}.} Let us take $\gamma=1$. Then, one of the solutions of~\eq{squareduffing} is
\[\xi(t) = y(t) = \delta \cn{\bigg[t\sqrt{1+\delta^2},\frac{\delta}{\sqrt{2(1+\delta^2)}}\bigg]}\, ,\]
where, as usual, $y(t)$ is the solution of \eq{ODE}-\eq{alphabeta} and it is explicitly given by \eq{soluzione_duffing}. Since this solution
is $T(\delta)$-periodic, its period is twice the period of the coefficient of the Hill equation \eq{squareduffing}. Thus, using the Floquet theory,
we deduce that the eigenvalues of the monodromy matrix for $\gamma=1$ are both $-1$. This shows that the line $\gamma=1$ is part of the boundary of $U_1$.\par
From \cite{abram} we recall that, for all $0<k<1$:
\[\frac{\partial \sn(u,k)}{\partial u} = \cn(u,k)\dn(u,k),\qquad
\frac{\partial \cn(u,k)}{\partial u} = -\sn(u,k)\dn(u,k),\]
\begin{equation}\frac{\partial \dn(u,k)}{\partial u} = -k^2\cn(u,k)\sn(u,k),\qquad\dn(u,k)^2-k^2\cn(u,k)^2=1-k^2.\label{idenjacobi}
\end{equation}
Let us take $k = \frac{\delta}{\sqrt{2(1+\delta^2)}}$ and, since there are no ambiguities, we put $\sn(u,k)=\sn(u)$.
Consider the function
\neweq{considerxi}
\xi(t) = \sn\Big(t\sqrt{1+\delta^2}\Big)
\endeq
so that, from the differential equalities \eq{idenjacobi}, we infer that:
\[\dot{\xi}(t) =\sqrt{1+\delta^2}\cn(t\sqrt{1+\delta^2})\dn(t\sqrt{1+\delta^2}),\]
\begin{eqnarray*}
	\ddot{\xi}(t) &=& -(1+\delta^2)\sn(t\sqrt{1+\delta^2})\Big(\dn(t\sqrt{1+\delta^2})^2+\frac{\delta^2}{2(1+\delta^2)}\cn(t\sqrt{1+\delta^2})^2\Big)\\
	&=& -(1+\delta^2)\sn(t\sqrt{1+\delta^2})\Big(\frac{2+\delta^2}{2(1+\delta^2)}+\frac{\delta^2}{1+\delta^2}\cn(t\sqrt{1+\delta^2})^2\Big)\\
	&=& -\bigg[1+\frac{\delta^2}{2} + \delta^2\cn\Big(t\sqrt{1+\delta^2}\Big)^2\bigg]\, \xi(t)\, .
\end{eqnarray*}
This shows that the choice \eq{considerxi} of $\xi(t)$ satisfies equation \eq{squareduffing} whenever $\gamma=1+\frac{\delta^2}{2}$.
Whence, on this parabola the eigenvalues of the monodromy matrix are $-1$. This proves that also the curve $\gamma=1+\delta^2/2$ is part of the
boundary of $U_1$.\par
Therefore, the boundary of $U_1$ consists of the line $\gamma=1$ and the parabola $\gamma=1+\delta^2/2$: the statement is so proved.\endproof

\noindent{\bf Proof of Theorem~\ref{striscia}.} In Theorem~\ref{exactsolutions} we showed that $\gamma=1$ is part of the boundary of $U_1$. Let $S_0$
be the first stability region for equation~\eqref{squareduffing}. Then, the line segment $\delta=0,\,0<\gamma<1$, is included in $S_0$.
Using~\cite[Chapter VIII, \S1-4]{yakubovich}, we conclude that $\gamma=1$ must be also part of the boundary of $S_0$.

From~\cite{abram} we recall that the following identity holds:
\begin{equation}\label{cos2sin2}
\cn(u,k)^2+\sn(u,k)^2=1\qquad(0<k<1).
\end{equation}
We take again $k=\frac{\delta}{\sqrt{2(1+\delta^2)}}$ and we drop it in the argument of $\sn$, $\cn$, $\dn$.
Then we consider the function
\neweq{considerxidn}
\xi(t) = \dn\Big(t\sqrt{1+\delta^2}\Big)\, .
\endeq
By arguing as for Theorem~\ref{exactsolutions} and by using~\eqref{idenjacobi} and~\eqref{cos2sin2}, we obtain:
\begin{equation*}
\dot{\xi}(t) = -\frac{\delta^2}{2\sqrt{1+\delta^2}}\cn(t\sqrt{1+\delta^2})\sn(t\sqrt{1+\delta^2})\, ,
\end{equation*}
\begin{align*}
\ddot{\xi}(t) &= -\frac{\delta^2}{2}\dn(t\sqrt{1+\delta^2})\big[-\sn(t\sqrt{1+\delta^2})^2 +\cn(t\sqrt{1+\delta^2})^2\big]\\
&=-\frac{\delta^2}{2}\big[-1+2\cn(t\sqrt{1+\delta^2})^2\big]\xi(t)\\
&=-\Big(-\frac{\delta^2}{2}+\delta^2\cn(t\sqrt{1+\delta^2})^2\Big)\xi(t)\, .
\end{align*}
This shows that~\eqref{considerxidn} satisfies equation~\eq{squareduffing} whenever $\gamma=-\frac{\delta^2}{2}$. Furthermore, this solution
has the same period $T(\delta)/2$ as the coefficient of the Hill equation~\eq{squareduffing}. Thus, using the Floquet theory, we conclude that
the characteristic multipliers of~\eqref{squareduffing} are both equal to $1$ when $\gamma = -\frac{\delta^2}{2}$.

From (4.2.i) p.60 in \cite{cesari}, we also know that if $p(t)\le0$, then the trivial solution of the Hill equation $\ddot{\xi}(t) +p(t)\xi(t)=0$ is
unstable. Therefore, if $\gamma\le-\delta^2$ the zero solution of~\eqref{squareduffing} is unstable. Thus, the parabola $\gamma=-\frac{\delta^2}{2}$ is
part of the boundary of $S_0$ and this concludes the proof. \endproof

\noindent{\bf Proof of Theorem~\ref{asymp}.} Assume that $(\delta,\gamma)\to(0,\ell^2)$. In order to apply the Burdina criterion we need asymptotic estimates for $\Phi(\delta,\gamma)$.
By neglecting $o(\delta^2)$, we obtain
\begin{eqnarray}
\Phi(\delta,\gamma) &=& 2\sqrt{2}\int_0^{\pi/2}\sqrt{\tfrac{\gamma+\delta^2\sin^2\theta}{2+\delta^2+\delta^2\sin^2\theta}}\, d\theta \notag \\
\ &\sim& 2\sqrt{\gamma}\int_0^{\pi/2}\left(1+\tfrac{\delta^2}{2\gamma}\sin^2\theta\right)
\left(1-\tfrac{\delta^2}{4}-\tfrac{\delta^2}{4}\sin^2\theta\right)\, d\theta \notag \\
\ &\sim& 2\sqrt{\gamma}\int_0^{\pi/2}\left(1-\tfrac{\delta^2}{4}+\tfrac{\delta^2}{2\gamma}\sin^2\theta-\tfrac{\delta^2}{4}\sin^2\theta\right)
\, d\theta \notag \\
\ &=& \pi\sqrt{\gamma}\left(1+\big(\tfrac{1}{\ell^2}-\tfrac32\big)\tfrac{\delta^2}{4}\right)\, .\label{stimona}
\end{eqnarray}

Next, we remark that
$$
\log\left(1+\frac{\delta^2}{\gamma}\right)\ \sim\ \frac{\delta^2}{\ell^2}\qquad\mbox{as }(\delta,\gamma)\to(0,\ell^2)\, .
$$
Combined with \eq{stimona}, this shows that the bounds in \eq{primavera} asymptotically become
$$
\gamma\le\ell^2+\left(\frac{3\ell^2}{4}-\frac12 -\frac{1}{\pi\ell}\right)\delta^2\ ,\quad
\gamma\ge\ell^2+\left(\frac{3\ell^2}{4}-\frac12 +\frac{1}{\pi\ell}\right)\delta^2\, .
$$
Whence, $U_\ell$ is asymptotically contained in the region where none of these two facts holds, that is, in the region defined by \eq{twoparabolas}.\par
Finally, the statement about the existence of $\delta_\gamma$ follows from what we have just proved and from Theorems~\ref{exactsolutions} and~\ref{striscia}.\endproof

\noindent{\bf Proof of Theorem~\ref{primalingua}.} If $\omega=1$, equations~\eqref{duffomega} and~\eqref{hill_scalata} are equivalent
to~\eqref{ODE} and~\eqref{squareduffing}, respectively. Therefore, it follows from the proof of Theorem~\ref{intornorigine} that
\[\xi(t) = \cn\bigg(t\sqrt{1+\delta^2},\frac{\delta}{\sqrt{2(1+\delta^2)}}\bigg)\]
is a solution of~\eqref{hill_scalata} whenever $\omega=1$. Thus the line $\omega=1$ is a resonant line and is part of the boundary of the resonance tongue
$U_1$. Let us now recall, from~\cite[Theorem 13]{bbfg}, that if $0<\omega<\big(\frac{21}{22}\big)^2$, then the trivial solution of~\eqref{hill_scalata} is
stable for all $\delta>0$. These two facts prove that the strip $(0,\infty)\times(0,1)$ is a stability region of the $(\delta,\omega)$-plane.
If $U_1=\{\omega:\tau(\delta)<\omega<\phi(\delta)\}$ then, by combining the previous observations with Proposition~\ref{known},
we infer that it has to be $\tau(\delta)\equiv1$, $\phi(\delta)>1$ for all $\delta$, and $\phi(\delta)\to3$ as $\delta\to\infty$.
The second inequality for $\phi(\delta)$ follows from Corollary~\ref{asympomega}.\endproof

\section{Remarks and open problems}

$\bullet$ In Figure \ref{monodromymnd1}, couples of resonant lines meet at the points $(\delta,\omega)=(0,\ell)$ with $\ell\in\N$. Between two
of these ``double points'' there is a gap of 1 in the $\omega$-direction. For $\delta=\infty$ the gap is larger and the $j$-th resonance line
emanates from the ``point'' $(\delta,\omega)=(\infty,\frac{j(j+1)}{2})$. Whence, below the line $\omega=s$ with $s\in\R_+\setminus\N$, there exist $2[s]$ resonant lines emanating from points on the $\omega$-axis and $[\frac{-1+\sqrt{1+8s}}{2}]$ resonant lines emanating from points at $\delta=\infty$. This
shows that there are ``more'' resonant lines for $\delta\to0$ than for $\delta\to\infty$. Here, $[\varrho]$ represents the integer part
of $\varrho$.\par\medskip\noindent
$\bullet$ Once a nontrivial solution $\xi_1$ of~\eqref{hill} is known, the standard method to find a second (linearly independent) solution $\xi_2$ yields
\neweq{xi2}
\xi_2(t)=\xi_1(t)\int^t\frac{d\tau}{\xi_1^2(\tau)}\, .
\endeq
Therefore, in the cases where we found an explicit solution of \eqref{hill}, see the proofs of Theorems
\ref{exactsolutions} and \ref{striscia}, one can also find another linearly independent explicit solution by using \eq{xi2} combined with
the computation of the indefinite integral of Jacobi elliptic function, see~\cite{carlson}. We omit the details.\par\medskip\noindent
$\bullet$ The plots in Figure \ref{monodromymnd1} suggest the following conjecture, which appears quite challenging to prove (or disprove): {\em all the
resonant lines emanating from $(0,\ell)$ (except $\omega\equiv1$) are graphs of strictly increasing functions with a unique flex point}.\par\medskip\noindent
$\bullet$ From \eq{twoparabolas} we infer that, when $\delta\to0$, the resonant tongues $U_\ell$ asymptotically lie between two upwards parabolas provided
that $\ell\ge2$. An interesting problem would be to investigate whether the gap between these two lines is of the order of $\delta^2$ or of higher
order $o(\delta^2)$ as $\delta\to0$.\par\medskip\noindent
$\bullet$ It would be interesting to investigate the stability for the following damped version of \eq{truebeam}:
$$u_{tt}+\rho u_t+u_{xxxx}-\tfrac{2}{\pi}\, \|u_x\|^2_{L^2(0,\pi)}\, u_{xx}=0\quad x\in(0,\pi)\, ,\ t>0\, ,$$
where $\rho>0$. How does $\rho$ modify the shape of the resonant tongues in Figure \ref{monodromymnd1}? We expect the tongues to retract in the
horizontal direction, but does a {\em uniform bound} $\eps_\rho>0$ exist such that all the nonlinear modes of \eq{truebeam} are linearly stable
whenever $\delta<\eps_\rho$? We point out that the simple argument of multiplying by an exponential, as for the Mathieu equation, does not work
for the Hill equation with squared Duffing coefficients.

\par\bigskip\noindent
\textbf{Acknowledgments.} This work is the continuation of the graduate dissertation of the first Author \cite{gaspa}.
The second Author is partially supported by the PRIN project {\em Equazioni alle derivate parziali di tipo ellittico
e parabolico: aspetti geometrici, disuguaglianze collegate, e applicazioni} and by the {\em Gruppo Nazionale per l'Analisi
Matematica, la Probabilit\`a e le loro Applicazioni (GNAMPA)} of the Istituto Nazionale di Alta Matematica (INdAM).

\end{document}